\documentclass[letterpaper]{article} 
\usepackage[accepted]{icml2019}  

\usepackage{hyperref}

\usepackage{lipsum}  

\usepackage{times}  
\usepackage{helvet}  
\usepackage{courier}  
\usepackage{url}  
\usepackage{graphicx}  
\usepackage{enumitem}
\usepackage{graphicx,psfrag,amsmath,amsfonts,amssymb, url, appendix}
\usepackage{algorithm, algorithmic}

\newtheorem{theorem}{Theorem}
\newtheorem{lemma}{Lemma}
\newtheorem{definition}{Definition}

\usepackage{booktabs} 
\usepackage{subcaption}

\def\reals{\mathbf{R}} 
\providecommand{\argmin}[1]{\mathop\mathrm{arg min}_{#1}}
\providecommand{\diag}{\mathop\mathrm{Diag}}

\def\norm#1{\left\|{#1}\right\|} 

\newcommand{\Id}{I}
\newcommand{\symmpd}[1]{{\mathbf{S}}^{#1}_{++}}  
\newcommand{\twonorm}[1]{\norm{#1}_2}
\newcommand{\onenorm}[1]{\norm{#1}_1}

\newcommand{\fronorm}[1]{\norm{#1}_{F}}

\newcommand{\minimize}[1]{\underset{#1}{\mathop{\rm minimize}} \quad ~}

\newcommand{\subjectto}{\mathop{\rm subject \ to \quad}}

\newcommand{\prox}{\mathbf{prox}}

\newcommand{\vertiii}[1]{{\left\vert\kern-0.25ex\left\vert\kern-0.25ex\left\vert #1
    \right\vert\kern-0.25ex\right\vert\kern-0.25ex\right\vert}}
 
\newcommand{\xhdr}[1]{\vspace{0.0mm}\noindent{\bf {#1}.\,}}

\usepackage{color}
\definecolor{Red}{rgb}{0.9,0.0,0.1}

\definecolor{redgray}{rgb}{0.5,0,0.40}

\definecolor{bluegray}{rgb}{0.15,0.20,0.40}
\definecolor{bluegraylight}{rgb}{0.35,0.40,0.60}
\definecolor{gray}{rgb}{0.3,0.3,0.3}
\definecolor{lightgray}{rgb}{0.7,0.7,0.7}
\definecolor{darkblue}{rgb}{0.2,0.2,1.0}
\definecolor{darkgreen}{rgb}{0.0,0.5,0.3}
\definecolor{myblue}{rgb}{0.627, 0.84, 0.97}

\begin{document}





\twocolumn[
\icmltitle{Variable Metric Proximal Gradient Method with
			Diagonal Barzilai-Borwein Stepsize}




\begin{icmlauthorlist}
\icmlauthor{Youngsuk Park}{stanford}
\icmlauthor{Sauptik Dhar}{lg}
\icmlauthor{Stephen Boyd}{stanford}
\icmlauthor{Mohak Shah}{lg}
\end{icmlauthorlist}

\icmlaffiliation{stanford}{Department of Electrical Engineering, Stanford University, Stanford, CA 94305.}
\icmlaffiliation{lg}{LG Sillicon Valley Lab, Santa Clara, CA 95050.}
\icmlcorrespondingauthor{Youngsuk Park}{youngsuk@stanford.edu}

\icmlkeywords{Machine Learning, ICML}

\vskip 0.3in
]



\printAffiliationsAndNotice{}  


\begin{abstract}

\textit{Variable metric proximal gradient} (VM-PG) is a widely used class of convex optimization method. Lately, there has been a lot of research on the theoretical guarantees of VM-PG with different metric selections. However, most such metric selections are dependent on (an expensive) Hessian, or limited to scalar stepsizes like the Barzilai-Borwein (BB) stepsize with lots of safeguarding. Instead, in this paper we propose an adaptive metric selection strategy called the \textit{diagonal Barzilai-Borwein (BB)  stepsize}. The proposed diagonal selection better captures the local geometry of the problem while keeping per-step computation cost similar to the scalar BB stepsize i.e. $O(n)$. Under this metric selection for VM-PG, the theoretical convergence is analyzed. Our empirical studies illustrate the improved convergence results under the proposed diagonal BB stepsize, specifically for ill-conditioned machine learning problems for both synthetic and real-world datasets.  

\end{abstract}

\section{Introduction}
\label{sec:intro}

We tackle a convex optimization in the composite form
\begin{align}
\minimize{x\in \reals^n} F(x):= f(x) + g(x), \label{cvxopt_composition}
\end{align}
where $x\in \reals^n$ is the decision variable, $f : \reals^n \rightarrow \reals$ is convex and differentiable, and $g : \reals^n \rightarrow \reals \cup \{\infty\}$ is convex and can be non-differentiable. Here, $g$ can be used to encode constraints on the variable $x$. Such structured form in \eqref{cvxopt_composition} appears across a wide range of machine learning problems like classification, regression, matrix completion etc. Proximal gradient methods have been widely adopted for solving the optimization problems involving such composite forms. There are several variants of the proximal gradient method in literature which offers various advantages such as improved computation costs, theoretical guarantees under mild conditions, practical rules for stepsize selections, etc \cite{tseng2000modified,combettes2005signal,barzilai1988two,zhou2006gradient,beck2009fast}. 
However, most of these modifications broadly follow a generic form known as \textit{Variable Metric Proximal Gradient method} (VM-PG) provided in Algorithm \eqref{algorithm:pg_vm} \cite{bonnans1995family,parente2008class}.
\begin{algorithm}
	\caption{Variable metric proximal gradient (VM-PG)}
	\label{algorithm:pg_vm}
\begin{algorithmic}
	\STATE {\bfseries given} a starting point $x^0 \in \reals^n$ 
   \REPEAT
   \STATE Update metric $U^k$
   \STATE $y^{k+1} = x^k - (U^k)^{-1} \nabla f(x^k)$
   \STATE $\begin{aligned} x^{k+1} &= \prox_{g,U^k}(y^{k+1})\\
    &:= \argmin{x} ~\left( g(x) + \frac12 \norm{y^{k+1}-x}_{U^k}^2 \right) 
    \end{aligned}$
   \UNTIL{stopping criterion $\twonorm{y^{k+1}- y^{k} }\leq \mathrm{\epsilon_{tol}}$ satisfied }
\end{algorithmic}
\end{algorithm}

Here, $x^k$ is the $k^{th}$ iterate, $U^k \in \symmpd{n}$ is a positive definite metric at the $k^{th}$ iteration, $\norm{z}_U =$ $ \sqrt{z^TUz}$ is the $U$-norm, and $\prox_{g,U}$ is the \textit{scaled} proximal mapping of $g$ relative to the metric $U$.


Note that, Algorithm \ref{algorithm:pg_vm} transforms to the standard proximal gradient algorithm for $U_k = (\alpha^k)^{-1} \Id$ where $\alpha^k$ is a scalar stepsize. And it becomes the proximal (quasi) Newton method for $U^k \approx \nabla^2 f(x^k)$ \cite{becker2012quasi,lee2014proximal}. These special cases have their respective pros and cons. For example, proximal Newton-type methods provide fast convergence in terms of iteration numbers but suffer worse per-step computation costs. On the other hand, proximal gradient methods have computationally attractive steps, but exhibit relatively slower convergence behaviors.

\xhdr{Summary of contributions} 
Even though many researchers have speculated that the usage of diagonal stepsizes would have superior convergence properties to that of scalar stepsize in general, few provide any practical diagonal stepsize rule across all convex optimization algorithms.  
In this paper, we propose a new adaptive rule for metric selection in VM-PG called \textit{diagonal Borzilai-Borwein} stepsize (Section 2). The proposed method tries to adopt the best of the two approaches: standard proximal gradient method and proximal Newton method. VM-PG with diagonal BB maintains low per-step computation cost $O(n)$ (similar to standard proximal gradient), while better satisfying a secant condition (i.e. better hessian approximation) at each iteration. This eventually leads to faster convergence behavior compared to standard proximal gradient methods (PG) with scalar BB stepsize. Convergence guarantees for the proposed method with line search is provided in Section 2. In Section 3, we introduce some computationally useful properties of the scaled proximal operator with (block) diagonal metric and derive the closed-form solutions for several interesting scaled proximal operations. Empirical results in Section 4 shows that the proposed VM-PG with diagonal metric provides better convergence than PG with the scalar BB stepsize. Conclusions are provided in Section 5.



\subsection{Related Work}

\xhdr{Spectral scalar stepsize} The BB Method \cite{barzilai1988two} is a popular approach for choosing a spectral stepsize in gradient descent methods for minimizing a quadratic objective. This method shows competitive convergence behavior compared to the widely used conjugate gradient method and demonstrates linear convergence behavior \cite{friedlander1998gradient}. In fact, this approach was later adopted for proximal gradient methods \cite{birgin2000nonmonotone,zhou2006gradient,wright2009sparse,goldstein2014field}. A special case of this method, also called as spectral projected gradient (SPG) or SpaRSA \cite{wright2009sparse}, demonstrates good numerical performance; even though its theoretical guarantees are not as strong as FISTA \cite{beck2009fast}. Recently, a similar idea was also proposed for the penalty parameters in Alternating Direction Method of Multiplier (ADMM) \cite{xu2016adaptive}. However, most of the these above approaches are limited to a scalar stepsize selection. Moreover, the empirical performance is often not in favor of non-quadratic problem, heavily depending on safeguarding parameters.

\xhdr{Variable non-scalar metric} The VM-PG (a.k.a variable metric forward-backward) method adopts a variable metric, rather than a scalar stepsize \cite{bonnans1995family}. This can provide better approximation of the local Hessian at each step $x^k$, which typically leads to improved convergence rates. There are several such metric rules proposed for both convex \cite{chouzenoux2014variable,salzo2016variable,lee2014proximal} and nonconvex  \cite{bonettini2016variable,boct2016inertial} problems. However, despite the theoretical convergence guarantees, most such proposed rules fall short in practical cases. For example, in the majorization-minorization principle \cite{chouzenoux2014variable,combettes2014forward}, deriving a majorization function compatible with the proximal operator is completely problem dependent. This hinders \textit{automatic} metric selection for many practical problems. Another example includes \cite{lee2014proximal}, where the metric updates approximate the Hessian similar to L-BFGS. This incurs an expensive Newton update followed by a scaled proximal step. Although, the Hessian approximations through rank $1$ updates on BB stepsize provide decent computational gains \cite{becker2012quasi}. However, scaled proximal mapping under this metric loses its closed-form solution property for many functions. And both gradient and proximal steps are not easily extendible to distributed algorithms. in brief, VM-PG incurs several computational limitations compared to standard proximal gradient, mainly for per-step computation costs.

\xhdr{Diagonal metric} A popular choice of a variable metric comes from the class of diagonal metrics. Such diagonal metrics are widely used in pre-conditioning strategies \cite{pock2011diagonal}. In fact, a specific form of diagonal metric has been succesfully applied for optimization over non-smooth objective functions in AdaGrad \cite{duchi2011adaptive}. Although widely used for a variety of problems, to our knowledge it has not been successfully applied to proximal gradient methods. In this paper we explore a diagonal (variable) metric for proximal gradient methods, and propose a new methodology for selecting the diagonal elements.

\section{Diagonal metric selection}
\label{sec:metric}


Our proposed adaptive rule for the spectral metric selection is motivated by the strengths and the shortcomings of the BB spectral stepsize method typically applied for gradient-type algorithms \cite{barzilai1988two}. First, we provide several insights into the BB method and highlight its limitations. Next, to alleviate the limitations, we propose the new adaptive diagonal metric selection strategy with convergence guarantees using a line search.


\subsection{Background and motivation} 

The proximal gradient step can be viewed as minimizing the overall function $F$ where the differentiable part $f$ is approximated into its second order form at $x^k$  (w.r.t some $ U^k \in \symmpd{n}$) \cite{chouzenoux2014variable}, 
\begin{align*}
&\prox_{g,U^k}(x^k - (U^k)^{-1} \nabla f(x^k) )= \\
&\argmin{x} ~g(x) + f(x^k) + \nabla f(x^k)^T(x - x^k)  + \frac{1}{2}\norm{x-x^k}_{U^k}^2.
\end{align*}
This motivates setting $U^k = \nabla^2 f(x^k)$ as a desirable choice following the proximal Newton method \cite{lee2014proximal}. However, using the Hessian typically incurs a high per-iteration cost. An alternative to that involves approximating the hessian using the \textit{secant condition},
\begin{align}
  U^k s^k \approx y^k \label{eq:secant_condition},
\end{align}
for the step $ s^k = x^{k} - x^{k-1}$ and the gradient change $y^k = \nabla f(x^{k}) - \nabla f(x^{k-1})$. 

\xhdr{Barzilai and Borwein (BB) method}
The \cite{barzilai1988two} BB method is a popular approach that estimates a scalar approximation of the Hessian by setting $U_k = (\alpha^k)^{-1} \Id$ which best satisfies \eqref{eq:secant_condition}. The two most widely used BB stepsizes are
\begin{align}
   & \alpha_{\mathrm{BB1}}^k := {\twonorm{s^k}^2 }/{\langle s^k, y^k \rangle}, \nonumber \\
   & \alpha_{\mathrm{BB2}}^k := {\langle s^k, y^k \rangle}/{\twonorm{y^k}^2 }\label{eq:bb12}.
\end{align}

\begin{definition} \label{def:smoothness}
A differentiable function $f : \mathbf{R}^n \rightarrow \mathbf{R}$ is $L$-smooth if
$\|\nabla f(x) -\nabla f(y)\|_2 \leq L \|x-y\|_2$ holds for all $x,y \in \reals^n$. 
And $f$ is $m$-strongly convex if $\langle \nabla f(x) -\nabla f(y), x-y \rangle \geq m \|x-y\|_2^2$ holds for all $x,y \in \reals^n$
\end{definition}

\begin{lemma} \label{bound_BB}
Let the differentiable $f$ be $L$-smooth and $m$-strongly convex. Then,  
\[
\frac{1}{L}\leq \alpha_{\mathrm{BB2}}^k\leq \alpha_{\mathrm{BB1}}^k\leq \frac{1}{m}.
\]
\end{lemma}
The proof follows from the definition of $L$-smoothness and $m$-strongly convexity and Cauchy-Schwartz inequality. Still for many degenerate scenarios with (large $L$ or small $m$) the bound is trivial and appropriate safeguarding for the numerical stability of the updates in eq. \eqref{eq:bb12} is still necessary. To this end, several modifications and safeguardings are adopted to the (original) BB stepsize \cite{zhou2006gradient,goldstein2014field}. One such numerical safeguarding on $s^k$ and $y^k$ is proposed using a hybrid choice between these two stepsizes following,
\begin{align}
\alpha_{\mathrm{BB}}^k&:=\alpha_{\mathrm{BB}}(s^k, y^k) \nonumber \\
	&=
	\begin{cases}
		\alpha_{\mathrm{BB2}}^k & \text{if }\alpha_{\mathrm{BB1}}^k < \delta\alpha_{\mathrm{BB2}}^k \\
			\alpha_{\mathrm{BB1}}^k - \frac{1}{\delta}\alpha_{\mathrm{BB2}}^k& \text{otherwise} 		
	\end{cases}, 	
	\label{eq:bb_hybrid}	
\end{align} 	
where the hyperparameter $\delta \in \reals$ is typically chosen as $2$. Lastly, if $\alpha_{\mathrm{BB}}^k$ in \eqref{eq:bb_hybrid} is negative, then the previous stepsize is selected, i.e.,  $\alpha_{\mathrm{BB}}^{k} = \alpha_{\mathrm{BB}}^{k-1}$. 

\xhdr{Caveats of scalar BB method}
Although, most such modifications and safeguardings are mainly designed to handle the instability in the (original) BB stepsize \eqref{eq:bb12} for ill-conditioned $f$. However, even with such modifications, the scalar BB may still be prone to inconsistencies. For example, note that $(\alpha_{\mathrm{BB1}}^k)^{-1}I$ and $(\alpha_{\mathrm{BB2}}^k)^{-1}I$ can be viewed as Hessian approximations in the Euclidean space. Under ill-conditioned settings, however, these \textit{scalar} approximations  may be far away from the true (non-Euclidean) Hessian geometry. Another case is that, after proximal mappings such as projections, the step ($s^k$) and gradient-change ($y^k$) directions can sometimes be close to being orthogonal. This causes degenerate scenarios with $\alpha_{\mathrm{BB1}}\rightarrow \infty$ or $\alpha_{\mathrm{BB2}} \rightarrow 0$. For such cases, the scalar estimates may significantly deviate from the secant condition \eqref{eq:secant_condition}, and in turn the Hessian geometry. 

\subsection{Diagonal Barzilai and Borwein stepsizes}
To better capture the Hessian geometry of $f$, we propose a diagonal metric $U^k$ at each iteration $k$ computed as follows
\begin{align} 
  &\minimize{u \in \reals^n} \norm{U s^k - y^k}_2^2 + \mu\fronorm{U - U^{k-1}}^2  \label{eq:dbb_opt}\\
  &\subjectto (\alpha_{\mathrm{BB1}}^k)^{-1} \Id \preceq U \preceq (\alpha_{\mathrm{BB2}}^k)^{-1} \Id, \nonumber \\
  & \qquad \qquad \quad ~~ U = \diag(u).  \nonumber
\end{align}
Here, the hyperparameter $\mu>0$ controls the trade-off between satisfying the secant condition \eqref{eq:secant_condition} and being consistent with the previous metric $U^{k-1}$. We choose a large $\mu$ if the Hessian does not change much over iterations. On the other hand, if Hessian changes fast, we choose a small $\mu$ which simply plays as a numerical safeguarding. Lastly, the diagonal elements are bounded by the (safeguarded) BB stepsizes in \eqref{eq:bb12}.

One advantage of the proposed formulation \eqref{eq:dbb_opt} is that it has a simple closed-form solution. For $U^k = \text{Diag}(u^k)$ and $u^k = [u_1^k,\ldots,u_n^k] \in \reals^n$, the solution to \eqref{eq:dbb_opt} is given as
\begin{align}
	&u^k_{i} = 
	\begin{cases}
		\frac{1}{\alpha_{\mathrm{BB1}}^k} \quad & \frac{s_{i}^ky_{i}^k + \mu u_{i}^{k-1} }{(s_{i}^k)^2 + \mu }  < \frac{1}{\alpha_{\mathrm{BB1}}^k}\\
		\frac{1}{\alpha_{\mathrm{BB2}}^k} & \frac{s_{i}^ky_{i}^k + \mu u_{i}^{k-1} }{(s_{i}^k)^2 + \mu } > \frac{1}{\alpha_{\mathrm{BB2}}^k}\\
		\frac{s_{i}^ky_{i}^k + \mu u_{i}^{k-1} }{(s_{i}^k)^2 + \mu } &\text{otherwise}
	\end{cases}\label{eq:dbb_solution}, 
\end{align}
where $s_i^k$ and $y^k_i$ are $i^{th}$ elements of $s^k$ and $y^k$ respectively.

\xhdr{Stability at degenerate scalar BB} 
The diagonal metric selection in \eqref{eq:dbb_solution} is likely to \textit{better satisfy} the secant condition \eqref{eq:secant_condition} compared to the (scalar) BB stepsize, whilst maintaining lower per-iteration cost compared to the proximal Newton-type methods (see Table \ref{table:cost}). For example, now in degenerate cases where $\langle s^k, y^k\rangle \approx 0$ (resulting $\alpha_{\mathrm{BB1}}^k\approx \infty$, $\alpha_{\mathrm{BB2}}^k\approx0$), the residual of secant condition with scalar BB \eqref{eq:bb_hybrid} can be very large. However, the residual $\|U^k s^k - y^k\|$ under diagonal metric can be much smaller for sufficiently small $\mu$. Moreover, $u^k$ at each iteration is still finite as long as $0 \leq u^{k-1} < \infty$ and $\mu > 0 $. This in practice, makes VM-PG with diagonal BB stepsize numerically more stable than the hybrid scalar BB in \eqref{eq:bb_hybrid}. 

In addition, although both (hybrid) scalar BB \eqref{eq:bb_hybrid} and diagonal BB depends on the previous metric, hybrid BB uses limited information wherein the previous value of the step-size is simply copied for the negative current stepsize. On the contrary, the diagonal BB better utilizes this additional information through a user-defined parameter $\mu$, casting an interplay between better Hessian approximations and/or numerical stability. For example,  setting large value of $\mu$ is the same as copying the previous step size (as adopted in hybrid scalar BB \eqref{eq:bb_hybrid}). The advantage of such a dynamic characterization of this interplay for different problem types is provided in the supplementary material.

\begin{table}[!t]
\centering
\begin{tabular}{|c|c|c|c|c|}
\hline
 & PG(BB) & \begin{tabular}[c]{@{}c@{}}VM-PG\\ (DBB)\end{tabular} & \begin{tabular}[c]{@{}c@{}}Prox \\ L-BFGS\end{tabular} & \begin{tabular}[c]{@{}c@{}}Prox \\ Newton\end{tabular} \\ \hline
Metric & $O(n)$ & $O(n)$ & $O(n^2)$ & $O(n^2)$ \\ \hline
Forward & $O(n)$ & $O(n)$ & $O(n^2)$ & $O(n^3)$ \\ \hline
\end{tabular}
\caption{Cost for computing metric $U^k$ and forward step ($ x^k - (U^k)^{-1} \nabla f(x^k)$).}
\label{table:cost}\vspace{-5pt}
\end{table}

\begin{figure}[t]
\centering
  \centering
  \includegraphics[width=0.8\linewidth]{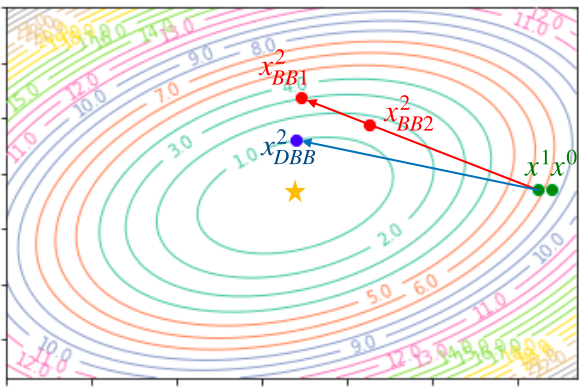}
  \caption{DBB (blue point) vs. BB 1 and BB 2 (red points): the three kinds of next iterate $x^2$ are pointed starting from initial $x^0$ and $x^1$ (green points)}
  \label{fig:vmpg_demo}
\end{figure}
As a simple illustration consider the toy example in Fig. \ref{fig:vmpg_demo}. Here, the magnitude of the diagonal BB stepsize (shown in blue) is bounded between BB 1 and BB 2 (shown in red); and the iterate direction using the diagonal approximation leads closer to the optimal solution. Hence, with similar per-step computation costs $O(n)$ (see \eqref{eq:dbb_solution}), the diagonal BB provides a better approximation of the Hessian and can eventually converge faster.


\xhdr{Remark} Note that, in essence a diagonal metric $U^k$ is equivalent to scaling the coordinates at each iteration, followed by a gradient and proximal step. Hence, the VM-PG (with diagonal metric) can be seen as performing a sequence of coordinate-scaling (or pre-conditioning); where the scale at each iteration $k$ depends on the local curvature (Hessian). This makes VM-PG with diagonal metric less sensitive to huge variations in the scale of the co-ordinates. 

Similar to many BB methods, VM-PG with the diagonal metric in \eqref{eq:dbb_solution} may still not guarantee convergence without line search (for penalized non-quadratic problems). Hence, we use \eqref{eq:dbb_solution} as an initial metric and additionally perform line search as shown in Algorithm \ref{algorithm:metric_selection}.



\subsection{Convergence under line search}

When $f$ is $L$-smooth, the standard proximal gradient method is guaranteed to converge for sufficiently small stepsize $\alpha < \frac{1}{L}$. Under no knowledge of the Lipschitz constant, there are several line search with backtracking \cite{boyd2004convex,beck2009fast} strategies which can still guarantee convergence. However, in practice a non-monotonic line search for VM-PG provides lower line search cost (per iteration), with better convergence results compared to the monotonic alternatives \cite{grippo1986nonmonotone,birgin2000nonmonotone,zhang2004nonmonotone,goldstein2014field}. 

\xhdr{Non-monotone line search}\label{sec:linesearch}
A non-monotonic line search allows the objective function $F(x)$ to increase between subsequent iterations, but results to an eventual decrease in its values. Here, given the current iterate $x^k$, an initial metric $U^k$ from \eqref{eq:dbb_solution}, and (a potential) next iterate $x^{k+1}$; the non-monotonic line search checks whether $(U^k, x^{k+1})$ satisfies the following criterion
\begin{align}
	F(x^{k+1}) \leq \hat F^k -\frac12 \norm{x^{k+1} - x^k}_{U^k}^2,\label{eq:linesearch}
\end{align}
where $M_{\mathrm{LS} }\geq 1$ is an integer line search parameter, and $\hat F^k= $ $\max\{F(x^k),F(x^{k-1})$ $,\ldots$ $,F(x^{k-\min(M_{LS}, k-1)}\})$. 
Then it backtracks by re-scaling the metric $U^k$ by a factor of $\beta >1$ until \eqref{eq:linesearch} is satisfied. 

\begin{algorithm}[!ht]
	\caption{VM-PG with diagonal BB metric}
	\label{algorithm:metric_selection}
\begin{algorithmic}
\STATE {\bfseries given} parameters $M_{\mathrm{LS} } \geq 1$, $\beta > 1$, $\mu>0$, a starting point $x^0, x^1\in \reals^n$, and initial metric $U^0\in \symmpd{n}$
\REPEAT
   \STATE Compute (safeguarded) $\alpha_{\mathrm{BB1} }^k$ and $\alpha_{\mathrm{BB2} }^k $ from  \eqref{eq:bb_hybrid}
   \STATE Initialize  $U^k$ from \eqref{eq:dbb_solution} 
   \STATE Update $x^{k+1} := \prox_{g,U_k}(x^k - (U^k)^{-1}\nabla f(x^k)) $
   \REPEAT
	   \STATE $U^k := \beta U^k$
	   \STATE $x^{k+1} :=\prox_{g,U_k}(x^k - (U^k)^{-1}\nabla f(x^k)) $
   \UNTIL{line search criterion in \eqref{eq:linesearch} is satisfied}\\
   \textbf{return} metric $U^k$ and next iterate $x^{k+1}$ 
\UNTIL{stopping criterion satisfied}
\end{algorithmic}
\end{algorithm}

 Next, we provide the convergence analysis for Algorithm \ref{algorithm:metric_selection}. We assume $f$ is $L$-smooth and $U^k >0$, then we have 
 

\begin{theorem} \label{thm:convergence} 
	For VM-PG in Algorithm \eqref{algorithm:metric_selection}, 
$F(x^k)$ converges to the optimal value $F^\star$, i.e., $\lim_{k\rightarrow \infty} F(x^k) := F^{\star}$. 
\end{theorem}
Additionally, the diagonal strategy in algorithm \eqref{algorithm:metric_selection} follows the following proposition,
\begin{theorem}
The VM-PG in Algorithm \eqref{algorithm:metric_selection}, with monotonic line search (i.e. $M_{\mathrm{LS}}=1$) satisfies,
	\begin{align*}
		\min_{k=1,\ldots, K}\norm{G_{U^k}(x^k)}^2_{(U^k)^{-1}} \leq \frac{2(F(x^0) - F^\star)}{K}
	\end{align*}
    where $G_{U^k}(x^k) \in \nabla f(x^k) + \partial g(x^k - (U^k)^{-1} \nabla f(x^k))$ and 
    $G_{U^k}(x^k)=0$ iff $0 \in \partial F(x^k)$.\\ \\
    Further, if $f$ is $m$-strongly convex, then 
    \begin{align*}
		\norm{x^{k+1} - x^\star}_{U^k}^2 
        \leq 
        (1-\frac{m}{u^k_{\mathrm{max}}})\norm{x^{k} - x^\star}_{U^k}^2 
	\end{align*}
    where, $u^k_{\mathrm{max}} = \max_i u_i^k$.
    
\end{theorem}
Proofs are provided in the supplementary material.







\section{Evaluation of scaled proximal mapping}
\label{sec:proximal_compute}


This section provides some useful properties of the scaled proximal mapping and illustrates the utility of such properties for machine learning algorithms.




\subsection{Properties of scaled proximal mapping} \label{sec:properties_proximal}
The key properties of the proximal mapping such as basic calculus, decomposition theorem, are  maintained for metric $U \in \symmpd{n}$.  


\begin{lemma}[Proximal Calculus \cite{rockafellar1976monotone}]
~
\begin{enumerate}
	\item If $f(x) = \alpha \phi(x)+b$, with $\alpha > 0$, then 
\begin{align*}
	\prox_{f, U}(x) = \prox_{\phi, U/\alpha}(x).
\end{align*}
	\item (Affine transformation) If $f(x) = \phi(Ax + b)$, with nonsingular $V \in \reals^{n \times n}$, then 
\begin{align*}
	\prox_{f, U}(x) = A^{-1} \left( \prox_{f, A^{-T}U A^{-1} }(Ax + b) - b\right).
\end{align*}
	\item \text{(Affine addition)} If $f(x) = \phi(x) + a^Tx + b$,
\begin{align*}
 	\prox_{f, U}(x) = \prox_{\phi, U}(x - U^{-1}a).
\end{align*}
	\item \text{(Regularization)} If $f(x) = \phi(x) + \frac12 \norm{x-a}_V^2$,
\begin{align*}
 	\prox_{f, U}(x) = \prox_{\phi, U+V}(x - (U+V)^{-1}(Ux+Va) ).
\end{align*}
	\item \text{Moreau decomposition \cite{becker2012quasi}}
\begin{align*}
x = \prox_{f, U}(x) + U^{-1} \prox_{f^*, U^{-1}}(Ux).
\end{align*}
\end{enumerate}
\end{lemma}

The next Lemma demonstrates the separability of proximal mapping for a separable function under block diagonal metric. This property enables distributing an algorithm (using consensus optimization).
\begin{lemma}[Separability]\label{lemma:separability}
Let $x = \{x_1,\ldots,x_{N}\}$ where $x_j \in \reals^{n_j}$, $U = \text{Blkdiag}(U_1,\ldots, U_N)$ where $U_j \in \symmpd{n_j}$, and $f$ be summable, meaning $f(x) = \sum_{j=1}^N f_j(x_{j})$.  Then the scaled proximal operator is separable, i.e., for each $j^{th}$ block,
\begin{align*}
\left( \prox_{f, U}(x) \right)_{j} = \prox_{f_j, U_j }(x_{j}).
\end{align*}
\end{lemma}

These properties provide practical utility for handling machine learning algorithms as illustrated next.  




\subsection{Examples of scaled proximal mapping}
Assume $\lambda, \lambda_1, \lambda_2 \in \reals_{+}$ are positive numbers, $U \in \symmpd{n}$, $u \in \reals^n_{+}$. We denote $(z)_i \in \reals$ as its $i$th element or $(z)_j \in \reals^{n_j}$ as its $j$th block under an explicit block structure, and $(z)_+ = \max(z,0)$. 

\xhdr{Lasso} For a lasso penalty $g(x) = \lambda \onenorm{x}$ and $U = \diag(u)$,
\begin{align*}
	\left( \prox_{g,U}(x) \right)_i = \text{sign}(x_i)(|x_i| - \lambda/u_i)_+.
\end{align*}

\xhdr{Group lasso} For a group lasso penalty $g(x) = \lambda \sum_{j}^N\twonorm{x_j}$ with $x_j \in \reals^{n_j}$ and $U = \text{Blkdiag}(u_1I_{n_1},\ldots,u_N I_{n_N})$,
\begin{align*}
	\left(\prox_{g,U}(x)\right)_j = \left(1 - \frac{\lambda}{u_j\twonorm{x_j}}\right)_+ x_j
\end{align*}

\xhdr{Elastic net} For a elastic net $g(x) = \lambda_1 \onenorm{x} + \lambda_2 \twonorm{x}^2$ and $U = \text{Diag}(u)$,
\begin{align*}
	\left( \prox_{g,U}(x) \right)_i = \text{sign}(x_i)\left(\frac{u_i}{\lambda_2 + u_i}|x_i| - \frac{\lambda_1}{\lambda_2 + u_i} )\right)_+ 
\end{align*}

\xhdr{Nonnegative constraint} Let $g(x) = \mathbf{1}(x\geq 0)$ be the nonnegative constraint. Then  
\begin{align*}
	\prox_{g,U}(x) = U^{-\frac12}(U^{\frac12}x)_+ 
\end{align*}
For $U = \text{Diag}(u)$, 
\begin{align*}
	\left( \prox_{g,U}(x)\right)_i = (x_i)_+ 
\end{align*}

\xhdr{Simplex constraint} Let $g(x) = \mathbf{1}(x\geq 0, \mathbf{1}^Tx=1)$ be the simplex constraint. Then 
for $U = \text{Diag}(u)$, 
\begin{align*}
	\left( \prox_{g,U}(x)\right)_i = (x_i - u_i^{-1}\nu)_+,
\end{align*}
Here, $\nu$ is the solution satisfying $\sum_i (x_i - u_i^{-1}\nu)_+ = 1$, which can be found efficiently via bisection on $\nu \in [\max_i u_i (y_i -1), \max_i u_iy_i]$.  

\xhdr{Consensus constraint} For $x = \{x_1,\ldots, x_{N_{\mathrm{node}}}\}$ with $x_j\in \reals^{n_j}$, let $g(x) = \delta_{\mathcal{C}}(x_1,\ldots, x_{N_{\mathrm{node}}})$ with a consensus constraint $\mathcal{C}=\{(x_1,\ldots, x_{N_{\mathrm{node}}}) \mid x_1=\ldots = x_{N_{\mathrm{node}}}\}$ and $\delta$ is a convex indicator. For $U = \text{Blkdiag}(U_1,\ldots, U_{N_{\mathrm{node}}})$ with $U_j \in \symmpd{n_j}$, 
\begin{align*}
	\left(\prox_{g,U}(x)\right)_j = \left(\sum_{j=1}^{N_{\mathrm{node}}} U_j\right)^{-1} \left(\sum_{j=1}^{N_{\mathrm{node}}} U_j x_j\right)
\end{align*}
where ${N_{\mathrm{node}}}$ is the number of nodes.

Note that all solutions can be computed with $O(n)$ cost. The derivations are provided in the supplementary material. 

\section{Experiments}
\label{sec:experiments}

This section provides the empirical results in favor of the proposed diagonal metric for VM-PG. We cover several applications with structure $F(x):=f(x) + g(x)$. 
\subsection{Applications}
%
\xhdr{Penalized quadratic programming} 
  \begin{align*}
      \minimize{x \in \reals^n}  \frac12 x^TQx + q^Tx + p + g(x),
  \end{align*}
where $Q \in \symmpd{n}$ and $q \in \reals^n$. For a regularizer, we use nonnegative constraint or lasso penalty with parameter $\lambda \in \reals_+$, i.e., $g(x) = \mathbf{1}_{\{z \mid z\geq 0\}}(x)$ or $g(x) = \lambda \onenorm{x}$ respectively. 
\\
\\
\xhdr{Penalized linear/logistic regression}
For $i=1 \ldots N$ samples of $a^{(i)} \in \reals^n$ and the associated label $b^{(i)}$, consider 
  \begin{align*}
      \minimize{x \in \reals^n}  \frac{1}{N}\sum_{i=1}^{N} l\left(x;a^{(i)}, b^{(i)}\right) +g(x),
  \end{align*}
where $l$ is a loss function, least square (linear) loss $l(\theta; a, b) = \twonorm{\theta^Ta - b}^2$ or logistic loss $l(\theta;a, b) = \log(1+e^{-b\theta^Ta}) $. Here, we also use nonnegative constraint or lasso penalty. 
\\

\subsection{Experimental setting}
We use several synthetic and real-world datasets. The experiment settings involve numerically challenging senarios with $N$ (sample size) $<<$ $n$ (feature dimension). A detailed discussion is provided next, 

\xhdr{Synthetic dataset} 

\begin{itemize} [leftmargin=*]\itemsep0em 
  \item For quadratic programming,  we consider well-conditioned ($\kappa = 10$) and ill-conditioned ($\kappa =10^4 $) cases. Here we use, $Q = HDH^T$ where $H$ is a random orthogonal matrix and $D = \diag(d_1,\ldots,d_n)$ with $\max_i d_i/\min_i d_i = \kappa$. 
 \item  For the penalized linear/logistic regression problems, we consider a small  ($N = 0.2n$) sample set generated from $a^{(i)} \sim \mathcal{N}(0, \Sigma)$ with some random $\Sigma \in \symmpd{n}$. Then the associated label $b^{(i)}$ is generated as follows. 
  \begin{itemize}[leftmargin=*]\itemsep0em 
    \item Least Square (LS) linear regression: $b^{(i)} = \left(a^{(i)}\right)^T x^\star + 0.2~v$ where $v \sim \mathcal{N}(0, I)$.
    \item Logistic regression (LR): $y = \sigma\left(\left(a^{(i)}\right)^T x^\star\right) + 0.2~w$ where $\sigma$ is sigmoid function $\sigma(z) = \log(1+e^{-z}) $ and $w \sim \text{Unif($0$, $1$)}$. Then take $b^{(i)}=1$ if $y\geq 0.5$ ~ or $b^{(i)}=-1$ otherwise. 
  \end{itemize}
\end{itemize}

\xhdr{Real-world datasets} 
We use two real-world datasets. Handwritten digit recognition MNIST 
\cite{lecun1998gradient} and object recognition CIFAR \cite{krizhevsky2012imagenet}. We show the results for a smaller subset of the dataset, illustrating the advantage of the proposed approach for highly ill-conditioned cases. The results using the entire dataset is provided in the supplementary material and show similar conclusions. For the MNIST and CIFAR datasets, we use LS and LR to estimate all labels (`0' - `9') and two labels (`1', `5') respectively.

\xhdr{Regularization parameter $\lambda$ and preconditioning}
For the synthetic and real-world datasets, we use $\lambda = 10^{-2}$ and $\lambda = 10^{-4}$ for LS and LR respectively. For the regression problems the data matrix $A \in \reals^{N \times n}$ is centered at $0$ and column-wise normalized to a unit $\ell_2$ norm.

\subsection{VM-PG algorithm parameters} 
Selecting the optimal $\mu$ (in eq. \eqref{eq:dbb_opt}) is problem dependent. For example of LR where the (local) Hessian may significantly change over iterations, a small $\mu$ allows the algorithm to properly capture the local geometry at each iteration and efficiently safeguard against degenerate cases. On the other hand, for cases like QP or LS where the local Hessian does not change over iterations, a large $\mu$ better captures the problem structure and is more desirable. Detailed experiments for different problem settings with varying $\mu$  are provided in the supplementary material. For the rest of this section we fix $\mu=10^{-6}$ to simplify our analysis. Also for the non-monotonic line search \eqref{eq:linesearch}, we set $M_{\mathrm{LS}} = 15$, $\beta = 2$, and adopt the modified stopping criterion with $\mathrm{\epsilon_{tol}} = 10^{-4}$ for QP/LS and $\mathrm{\epsilon_{tol}}=10^{-2}$ for LR problems following \cite{goldstein2014field}.

\begin{table*}[!h]
\centering
\scalebox{0.9}{

\begin{tabular}{|c|c|c|c|c|c|c|c|}
\hline
$f(x)$ & QP & QP & $f(x)$ & LS & LR & LS & LR \\
$g(x)$ & nonneg. & nonneg. & $g(x)$ & nonneg. & nonneg. & lasso & lasso \\
($\kappa$, $n$) & ($10$, $10^3$) & ($10^5$, $10^3$) & ($N$, $n$) & ($200$, $10^3$) & ($200$, $10^3$) & ($200$, $10^3$) & ($200$, $10^3$) \\ \hline \hline
PG(BB) & 9.8 (0.3) & 22.1 (0.62) & PG (BB) & 52.3 (1.22) & 54.5 (1.71) & 82.1 (2.09) & 61.5 (2.24) \\ \hline
\begin{tabular}[c]{@{}c@{}}VM-PG \\ (Diagonal BB)\end{tabular} & 8.2 (0.27) & \bf{16.2 (0.49)} & \begin{tabular}[c]{@{}c@{}}VM-PG\\ (Diagonal BB)\end{tabular} & 46.15  (1.08) & \bf{46.2 (1.27)} & 84.9 (2.21) & \bf{45.5 (1.13)} \\ \hline
\end{tabular}

}
\caption{Average number of iterations (CPU times in sec) for the convergence of VM-PG (Diagonal BB) and PG (BB) for penalized quadratic programming (QP), least square (LS), and logistic regression (LR) with synthetic dataset.}
\label{table:vmpg}\vspace{-5pt}
\end{table*}

\begin{table*}[!h]
\centering
\scalebox{0.9}{

\begin{tabular}{|c|c|c|c|c|}
\hline
$f(x)$ & LS & LR & LS & LR \\
$g(x)$ & lasso & lasso & lasso & lasso \\
\hline \hline
\begin{tabular}[c]{@{}c@{}}Data\\ ($N$, $n$)\end{tabular} & \begin{tabular}[c]{@{}c@{}}MNIST\\ ($240$, $784$)\end{tabular} & \begin{tabular}[c]{@{}c@{}}MNIST\\ ($1250$, $784$)\end{tabular} & \begin{tabular}[c]{@{}c@{}}CIFAR\\ ($625$, $3072$)\end{tabular} & \begin{tabular}[c]{@{}c@{}}CIFAR\\ ($500$, $3072$)\end{tabular} \\ \hline
PG (BB) & 83 (2.24) & 181 (5.52) & 175 (5.7) & 91 (4.42) \\ \hline
\begin{tabular}[c]{@{}c@{}}VM-PG\\ (Diagonal BB)\end{tabular} & 78 (2.01) & \bf{133 (3.83)} & 181 (5.52) & \bf{49(2.67)} \\ \hline
\end{tabular}

}
\caption{Iterations (CPU times in sec) for the convergence of VM-PG (Diagonal BB) and PG (BB) for $\ell_1$ penalized least square (LS), and logistic regression (LR) with subsampled real datasets. }
\label{table:vmpg_real}\vspace{-5pt}
\end{table*}

\begin{figure*}[!h]
\centering
\begin{subfigure}[t]{.3\textwidth}
  \centering
  \includegraphics[width=1\linewidth]{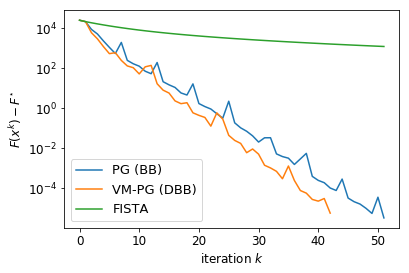}
  \caption{Nonneg. LS for synthetic data with $N =200,n = 1000.$ }
  \label{fig:sub1}
\end{subfigure} \quad
\begin{subfigure}[t]{.3\textwidth}
  \centering
  \includegraphics[width=1\linewidth]{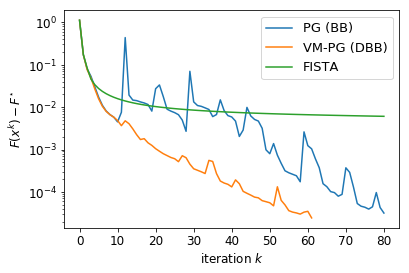}
  \caption{$\ell_1$ penalized LR for synthetic data with $N =200,n = 1000$.}
  \label{fig:sub2}
\end{subfigure} \quad
\begin{subfigure}[t]{.3\textwidth}
  \centering
  \includegraphics[width=0.9\linewidth]{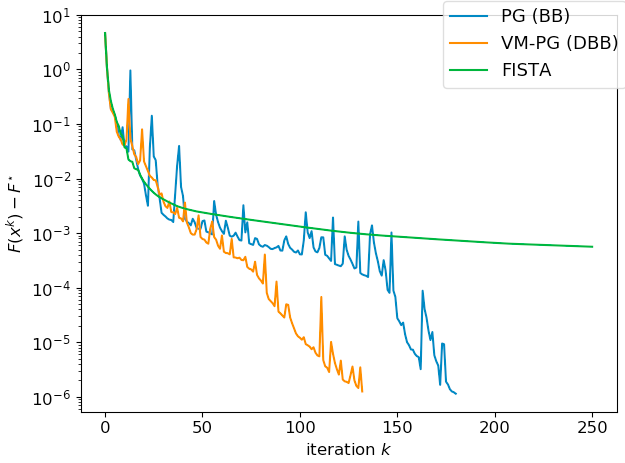}
  \caption{$\ell_1$ penalized LR for MNIST data with $N = 1250 , n =784$.}
  \label{fig:sub3}
\end{subfigure}\\


\begin{subfigure}[t]{.3\textwidth}
  \centering
  \includegraphics[width=1\linewidth]{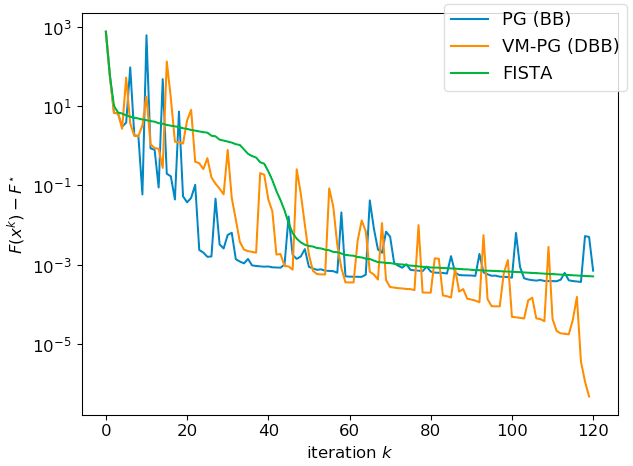}
  \caption{$\ell_1$ penalized LS for MNIST $N = 240,n =784$}
  \label{fig:sub4}
\end{subfigure}\quad%
%
\begin{subfigure}[t]{.3\textwidth}
  \centering
  \includegraphics[width=1\linewidth]{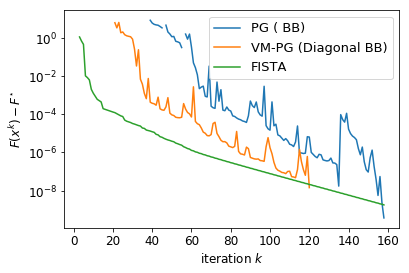}
  \caption{$\ell_1$ penalized LR for CIFAR with $N=10000,n =3072$}
  \label{fig:sub5}
\end{subfigure}\quad
\begin{subfigure}[t]{.3\textwidth}
  \centering
  \includegraphics[width=1\linewidth]{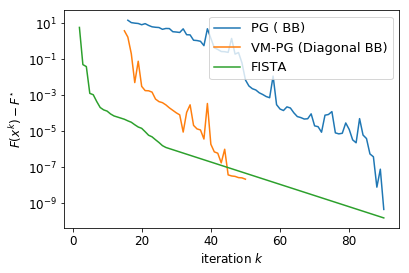}
  \caption{$\ell_1$ penalized LR for CIFAR with $N =500,n =3072$}
  \label{fig:sub6}
\end{subfigure}\\
%
%

\vspace{0mm}
\caption{Typical convergence behaviours of VM-PG with diagonal BB stepsize (orange), PG with BB stepsize (blue), and Accelerated PG (FISTA) (green): (a) and (b) are for synthetic data, (c) and (d) are for the MNIST dataset, and (e) and (f) are for the CIFAR dataset}
\label{fig:test}
\end{figure*}

\subsection{Results} 
Table \ref{table:vmpg} shows the total number of iterations (and CPU times in sec) for the convergence of the VM-PG (with diagonal BB) vs. PG (scalar BB), averaged over 100 experimental runs for the synthetic data. And Table \ref{table:vmpg} show the results for subsampled (ill-conditioned) MNIST and CIFAR dataset. The results the accelerated proximal gradient method (FISTA) (a non BB-type method) \cite{beck2009fast}, is provided only as a reference. Note that, all the three methods require similar per iteration computational costs, i.e., $O(n)$ to compute and store metric, $O(mn)$ or $O(n^2)$ cost for gradient step, $O(n)$ for proximal step; of which the gradient steps are dominant. Hence, the Fig. \ref{fig:test} majorly illustrate the convergence behaviors of the methods in terms of their iteration counts.


\xhdr{Penalized QP} As seen from Table \ref{table:vmpg}, for well-conditioned $Q$ with $\kappa \sim 10$, both the methods provide fast convergence without any significant difference. However, VM-PG with the diagonal BB selection \eqref{eq:dbb_solution} significantly outperforms standard PG (BB) (with $\sim 20\% $ computation improvement) for the ill-conditioned $Q$ ($\kappa \sim 10^4$). For unconstrained QP too we see similar results. In fact, VM-PG exhibits less oscillation and requires lower line search iterations. Additional figures illustrating such convergence properties are provided in the supplementary material.

\xhdr{Penalized regression}
As seen from Table \ref{table:vmpg} the VM-PG (diagonal BB) significantly outperforms the standard PG (scalar BB). Additional convergence behavior for both real/synthetic datasets are provided in Fig. \ref{fig:test}. Fig. \ref{fig:test} shows that the proposed VM-PG (with diagonal BB) significantly outperforms PG (scalar BB) with $\sim 20 \%$ improvement for LR lasso. In fact, the proposed VM-PG (diagonal BB) performs as good as (or even better) than FISTA in most of the cases (see \ref{fig:test}). For cases where FISTA outperforms the VM-PG (diagonal BB), optimally tuning the  $M_{\mathrm{LS}}$ parameter provides significant improvement for VM-PG (diagonal BB). A more detailed study on comparisons with other state-of-art methods like FISTA under different problem settings and the equivalent parameter optimizations for VM-PG (diagonal BB) is an open research problem. Additional results including experiments on entire dataset confirm the results presented in this section and are provided in the supplementary materials.


In short, the results illustrate that for ill-conditioned problems, the proposed VM-PG with diagonal BB better captures the local geometry of the problem and leads to better convergence results, compared to PG with scalar BB.

\section{Conclusion}
\label{sec:conclusion}
This paper proposes a diagonal BB metric for the variable proximal gradient method. The proposed diagonal metric provides a better estimate of the ill-conditioned local Hessian compared to the standard scalar BB approach, resulting to a faster convergence. Combined with a nonmonotonic line-search the overall algorithm is guaranteed to converge. Finally, for several machine learning applications with synthetic and real-world datasets, empirical results exhibit improved convergence behavior for the proposed methodology.



\nocite{langley00}

\bibliography{refs}
\bibliographystyle{icml2019}



\appendix
{\centering
\section*{Appendix}
}
\label{sec:supplementary}

\section{Proofs and Derivations}

\subsection{Proofs for Theorems 1 \& 2}
For the proofs of the theorems 1 and 2 we first provide the following lemmas, 
\begin{lemma}\label{lemma:proximal} For any proximal mapping the following holds,
\[
y = \prox_{g,U}(x) \text{ if and only if }  U(x-y)  \in \partial g(y)
\]
\end{lemma}
\textbf{Proof}: For any $y$ that minimizes $\prox_{g,U}(x)$ we have,  
\begin{flalign}
&y = \underset{v}{\text{argmin}} \; g(v)+(1/2)||x-v||_U^2  \quad (\text{from definition}) \nonumber && \\
&\Leftrightarrow 0 \in \partial g(y) + U(x-y) \quad (\text{at minima})&& \nonumber
\end{flalign} 

\begin{lemma}\label{lemma:key_inequality}
For any $x^+ = \prox_{g,U}(x - U^{-1} \nabla f(x) ) = x - U^{-1} G_U(x)$  \; where \; $G_U(x) = U(x - \prox_{g,U}(x - U^{-1} \nabla f(x) ) )$ and $U \succeq L \cdot \Id$ we have $\forall z \in \reals^n$,
\begin{align*}
F(x^+) 
&\leq  F(z) + G_U(x)^T(x-z) -\frac{m}{2}\norm{z-x}_2^2 
\\
& \qquad - \frac{1}{2} \norm{ G_U(x)}_{U^{-1} }^2 \quad \quad 
\end{align*}
\end{lemma}
\textbf{Proof:}
For $x^+  = x - U^{-1} G_U(x) $
we have
\begin{align*}
g(x &- U^{-1} G_U(x) ) \\
&\overset{(a)}{\leq} g(z) - \partial g(x - U^{-1} G_U(x) )^T(z - x + U^{-1}G_U(x))\\
&\overset{(b)}{=}  g(z) -(G_U(x) - \nabla f(x))^T(z - x + U^{-1}G_U(x))\\
&\overset{}{=}  g(z) + G_U(x)^T(x-z) - \norm{G_U(x)}_{U^{-1}}^2 \\
&\qquad + \nabla f(x)^T(z - x + U^{-1}G_U(x))
\end{align*}
where (a) holds due to convexity of $g$, (b) holds from Lemma \ref{lemma:proximal}. 
Next, 
\begin{align*}
f(x &- U^{-1} G_U(x) ) 
\\
&\overset{(c)}{\leq } f(x) - \nabla f(x)^T U^{-1} G_U(x) + \frac{L}{2} \norm{U^{-1} G_U(x)}_2^2 \\ 
&\overset{(d)}{\leq } f(x) - \nabla f(x)^T U^{-1} G_U(x) +\frac{1}{2} \norm{ G_U(x)}_{U^{-1} }^2\\ 
&\overset{(e)}{\leq } f(z) -  \nabla f(x)^T(z-x ) -\frac{m}{2}\norm{z-x}_2^2\\
&\quad  - \nabla f(x)^T U^{-1} G_U(x) + \frac{1}{2} \norm{ G_U(x)}_{U^{-1} }^2\\ 
&\overset{}{=} f(z) -  \nabla f(x)^T(z-x + U^{-1} G_U(x)) -\frac{m}{2}\norm{z-x}_2^2 \\
&\qquad + \frac{1}{2} \norm{ G_U(x)}_{U^{-1} }^2\\ 
\end{align*}
where (c) holds due to $L$-smoothness, (d) holds by $U \succeq L \cdot \Id$, (e) holds due to $m$-strongly convexity.

Therefore,
\begin{flalign}
F(x^+) & = f(x - U^{-1} G_U(x) ) + g(x - U^{-1} G_U(x) ) \nonumber &&\\
&\leq  F(z) + G_U(x)^T(x-z) -\frac{m}{2}\norm{z-x}_2^2 - \frac{1}{2} \norm{ G_U(x)}_{U^{-1} }^2 \; \nonumber && 
\end{flalign}  

\begin{lemma}\label{lemma:descent_lemma}
For the updates in Algorithms 1 \& 2 where, $x^{k+1} = prox_{g,U^k}(x_k-(U^k)^{-1} \nabla f(x^k))$ assuming $U^k \succeq L \cdot I$ we have $\forall k$,
\begin{align*}
F(x^{k+1}) \leq  F(x^k)  - \frac{1}{2} \norm{ G_{U^k}(x)}_{(U^k)^{-1} }^2
\\
=  F(x^k)  - \frac{1}{2} \norm{ x^{k+1} - x^k}_{U^k }^2
\end{align*}
\end{lemma}
\xhdr{Proof:} The proof follows by setting $z = x^k$, $x = x^k$ and $x^+ = x^{k+1}$ in Lemma \ref{lemma:key_inequality}.

\begin{lemma} \label{lemma:back_track}
	Assuming $f$ is $L$-smooth the linesearch criterion (see eq. (7)) in Algorithm 2 is satisfied within finite number of backtrackings.
\end{lemma} 

\textbf{Proof:}
From algorithm 2 we have $U^k>0$. Hence with finite number of backtracking using $\beta > 0$ we can have $U^k \succeq L \cdot I$. This ensures, 
\begin{flalign}
&F(x^{k+1}) \leq F(x^k) - (1/2) \|x^{k+1} -x^k \|_{U^k} \quad (\text{Lemma \ref{lemma:descent_lemma}}) && \nonumber \\
&\Rightarrow F(x^{k+1}) \leq \hat F^k - (1/2) \|x^{k+1} -x^k \|_{U^k} \quad (\because F^k \geq F(x^k))&& \nonumber
\end{flalign} 

With the above Lemmas in place we prove the main Theorems 1 \& 2. For readability we re-write the theorems here,

\begin{theorem}
For VM-PG in Algorithm 2, 
$F(x^k)$ converges to the optimal value $F^\star$, i.e., $\lim_{k\rightarrow \infty} F(x^k) := F^{\star}$.
\end{theorem}

\textbf{Proof:} 
Since the iterates satisfy the linesearch criteria, following Lemma \ref{lemma:descent_lemma} and \ref{lemma:back_track} we have,
\begin{align} 
	F(x^{k+1}) \leq \hat F^k - \frac12 \norm{x^{k+1} - x^k}_{U^k}^2 \nonumber
    \\
     = \hat F^k - \frac{1}{2} \norm{ G_{U^k}(x)}_{{(U^k)}^{-1} }^2 \nonumber
\end{align}

Further, $\{\hat F^k\}$ is monotonically decreasing sub-sequence of $\{F(x^{k})\}$. Let this sub-sequence be indexed as $\{F(x^{k'(i)})\}$ for some $k - M_{\mathrm{LS}}\leq k'(i)\leq k$. For this sub-sequence at limit we have, 
\begin{flalign} \label{eq_limitpoint}
0 &= 
\lim_i \norm{x^{k'(i)+1} - x^{k'(i)}}_{U^{k'(i)}}^2  && \nonumber \\
&=
\lim_i \norm{ G_{U^{k'(i)}}(x^{k'(i)})}_{{(U^{k'(i)})}^{-1} }^2 && 
\end{flalign}
Also, assuming that the limit point exist, let this limit point be $ \hat F^\star = F(\hat x^\star)$. For this limit point, \eqref{eq_limitpoint} implies $G_{U^\star}(\hat x^\star) = 0$ ($\because U^k >0; \forall k$). But we know, $G_{U^\star}(\hat x^\star) = 0$ iff $0 \in \partial F(\hat x^\star)$. Hence, this limit point $\hat x^\star$ is a stationary point of $F(x)$. Finally, $\hat x^\star$ is also the global minima under convexity of $F$. 

\begin{theorem} 
The VM-PG in Algorithm 2, with monotonic line search (i.e. $M_{\mathrm{LS}}=1$) satisfies,
	\begin{align*}
		\min_{k=1,\ldots, K}\norm{G_{U^k}(x^k)}^2_{(U^k)^{-1}} \leq \frac{2(F(x^0) - F^\star)}{K}
	\end{align*}
    where $G_{U^k}(x^k) \in \nabla f(x^k) + \partial g(x^k - (U^k)^{-1} \nabla f(x^k))$ and 
    $G_{U^k}(x^k)=0$ iff $0 \in \partial F(x^k)$.\\ \\
    Further, if $f$ is $m$-strongly convex, then 
    \begin{align*}
		\norm{x^{k+1} - x^\star}_{U^k}^2 
        \leq 
        (1-\frac{m}{u^k_{\mathrm{max}}})\norm{x^{k} - x^\star}_{U^k}^2 
	\end{align*}
    where, $u^k_{\mathrm{max}} = \max_i u_i^k$.
\end{theorem}

\textbf{Proof:}
For the first part, from Lemma \ref{lemma:descent_lemma}, 
\begin{align*}
F(x^{k+1}) \leq  F(x^k)  - \frac{1}{2} \norm{ G_U(x^k)}_{(U^k)^{-1} }^2 \; \forall k
\end{align*}
Reordering terms and averaging over iterations $k = 1 \ldots K$ gives,
\begin{align*}
\frac{1}{K} \sum_{k=1}^K\norm{ G_{U^k}(x^k)}_{(U^k)^{-1} }^2 
&\leq \frac{2}{K} \sum_{k=1}^K F(x^k)  - F(x^{k+1}) 
\\
&\leq \frac{2(F(x^0)  - F(x^\star)) }{K} .
\end{align*}
And LHS is lower bounded by 
\[
\frac{1}{K} \sum_{k=1}^K\norm{ G_{U^k}(x^k)}_{(U^k)^{-1} }^2 
\geq
\underset{{k=1,\ldots,K}}{\min }\norm{ G_{U^k}(x^k)}_{(U^k)^{-1} }^2 .
\]

For the second part, substituting $z = x^\star$ in Lemma \ref{lemma:key_inequality} gives,
\begin{align*}
F&(x^+) - F^\star  \leq G_U(x)^T(x-x^\star) -\frac{m}{2}\norm{x-x^\star}_2^2 - \frac{1}{2} \norm{ G_U(x)}_{U^{-1} }^2 \\
& = \frac12 \left(\norm{x - x^\star}_U^2  - \norm{x - x^\star - U^{-1}G_U(x) }_U^2  -m\norm{x-x^\star}_2^2\right)  \\
& = \frac12 \left(\norm{x - x^\star}_U^2  - \norm{x^+ - x^\star }_U^2  - m \norm{x-x^\star}_2^2 \right).
\end{align*}
Reordering terms give
\begin{align*}
\norm{x^+ - x^\star }_U^2 
&\leq \norm{x - x^\star}_U^2  - \left( 2 (F(x^+) - F^\star) + m \norm{x-x^*}_2^2 \right) \\
&\overset{}{\leq } (1 - \frac{m}{u_{\mathrm{max}}^k})\norm{x - x^\star}_U^2  
\end{align*}
where last inequality holds due to $F(x^+) - F^\star \geq 0$ and $\norm{a}_2^2 \geq (1/u_{\mathrm{max}}^k)\norm{a}_U^2$
for $u_{\mathrm{max}}^k = \max_{i=1,\ldots,n} u_i^k$. 

\subsection{Derivations for the proximal forms of the constraints in Section 3.2}
Most of the derivations are immediate from the properties in Section 3.1. Here, we derive the non-trivial case of simplex constraint. \\
\xhdr{Simplex constraint}
By taking dual, 
\begin{align*}
    \max_{\lambda \geq 0, \nu} \min_x \frac12 \|x-y\|_U^2  + \lambda ^T(-x) + \nu^T(\mathbf{1}^Tx-1)
\end{align*}
Taking derivative gives 
\[
U(x-y) - \lambda + \nu \mathbf{1}=0
\]
and then the dual becomes 
\begin{align} \label{eq_dual_simplex}
    \max_{\lambda \geq 0, \nu} &\frac12 \|\lambda - \nu \mathbf{1} \|_U^2  + (y + U^{-1}(\lambda - \nu \mathbf{1}))^T (-\lambda + \nu \mathbf{1}) - \nu
    \nonumber \\
    =
    \max_{\lambda \geq 0, \nu}& -\frac12 \|U y + (\lambda - \nu \mathbf{1}) \|_{U^{-1}}^2   - \nu
\end{align}
The optimal solution for \eqref{eq_dual_simplex} is \begin{flalign} 
\lambda_i = \left\{
\begin{array}{l l}
    -(u_iy_i - \nu);\quad  (u_iy_i - \nu)<0 \\
    0;\quad  \text{else}
\end{array}\right.  \nonumber \end{flalign} This gives the dual 

\begin{align*}
    \max_{\nu} -\frac12 \|[U y + - \nu \mathbf{1}]_+ \|_{U^{-1}}^2   - \nu
\end{align*}
Taking derivative over $\mu$ gives $
\sum_i (y_i - u_i^{-1}nu)_+ =1
$. Thus we get the solution $x_i = (y_i - u_i^{-1}\nu)_+$.  

\section{Additional Results} \label{sup:algo_param}

\subsection{Convergence behaviour of the VM-PG algorithm with varying $\mu$ values} 

\begin{figure}[H]
\centering
\begin{subfigure}{.4\textwidth}
  \centering
  \includegraphics[width=1\linewidth]{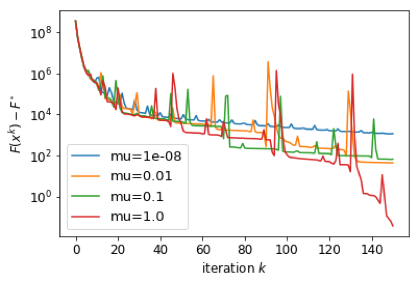}
  \caption{No linesearch}
  \label{fig:sub11}
\end{subfigure}
\begin{subfigure}{.4\textwidth}
  \centering
  \includegraphics[width=1\linewidth]{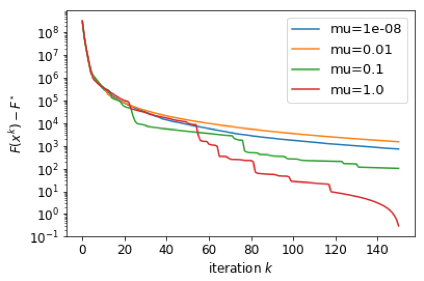}
  \caption{Monotonic linesearch}
  \label{fig:sub12}
\end{subfigure}
\caption{The effect varying the $\mu$ parameter in VM-PG algorithm for unconstrained QP problem with $\kappa = 10^5$.}
\vspace{-5pt}
\label{fig:test_Quad}
\end{figure}

For a quadratic programming problem as discussed in Section 4.1 without the constraints. The typical convergence behaviour for varying $\mu = [1e^{-8}, 0.01, 0.1 ,1]$ values are shown in Fig. \ref{fig:test_Quad}. The results are generated using similar experimental settings (with $\kappa = 10^4$) discussed in section 4.2. Fig. \ref{fig:test_Quad}(a) shows the convergence behaviour without any linesearch.  Fig. \ref{fig:test_Quad} (b) provides the results using monotonic linesearch. As seen from the figures, a larger $\mu$ value illustrates improved convergence behaviour. This can be attributed to the fact that the local Hessian for the  QP problem does not change significantly over iterations. Hence, putting higher weight on the second term in eq. 2 ensures this and provides improved convergence results.    

\begin{figure}[h]
\centering
\begin{subfigure}{.4\textwidth}
  \centering
  \includegraphics[width=1\linewidth]{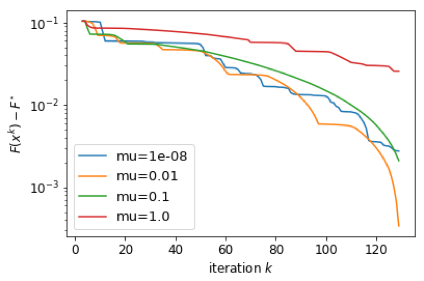}
  \caption{$\ell_1$ penalized LR with monotonic line search}
  \label{fig:sub13}
\end{subfigure}
\begin{subfigure}{.4\textwidth}
  \centering
  \includegraphics[width=1\linewidth]{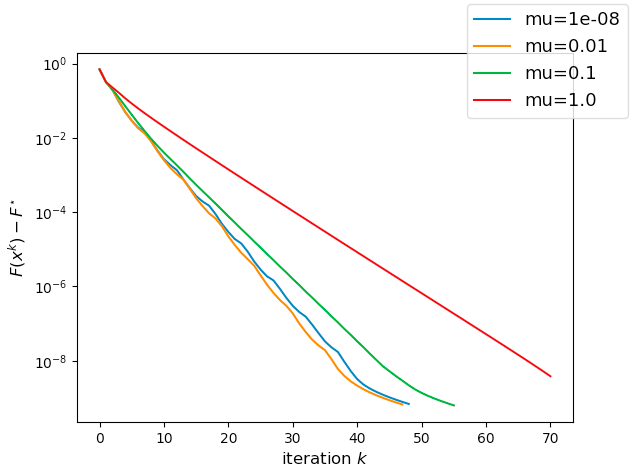}
  \caption{nonnegative LR with monotonic line search}
  \label{fig:sub14}
\end{subfigure}\\
\vspace{0mm}
\caption{The effect varying the $\mu$ parameter in VM-PG algorithm for logistic regression problems. }
\vspace{-5pt}
\label{fig:pen_logreg}
\end{figure}

Next we provide an analysis for varying $\mu$ for the penalized logistic regression (LR) problem in section 4.1. Fig \ref{fig:pen_logreg} shows the convergence behaviour using similar experimental settings (with $N = 200, n=10^3$) as discussed in section 4.2. As seen from the figures, the performance using smaller $\mu \leq 0.01$ values is better than that using larger $\mu \geq 1$ values. For problems like LR where the (local) Hessian may change abruptly over iterations, a small $\mu \leq 10^{-8}$ is preferable (also confirmed from the results). Such a selection ensures numerical safeguarding and does not heavily depend heavily on the previously estimated metric. 

In short, the $\mu$ parameter enables us to incorporate the information about the problem structure and hence improve the overall performance of the algorithm.



\begin{figure*}
\centering
\begin{subfigure}{.4\textwidth}
  \centering
  \includegraphics[width=1\linewidth]{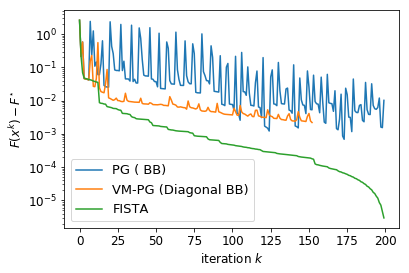}
  \caption{$\ell_1$ penalized LR for MNIST using
      ($N, n$) = ($12000, 784$)}
  \label{fig:sub4}
\end{subfigure} \hspace{0.05\textwidth} %
\begin{subfigure}{.4\textwidth}
  \centering
  \includegraphics[width=1\linewidth]{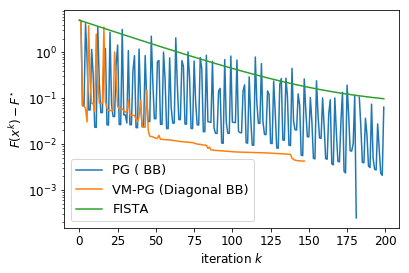}
  \caption{$\ell_1$ penalized LS for MNIST using
  		($N,n$)=($60000,784$)}
  \label{fig:sub6}
\end{subfigure}\\
\begin{subfigure}{.4\textwidth}
  \centering
  \includegraphics[width=1\linewidth]{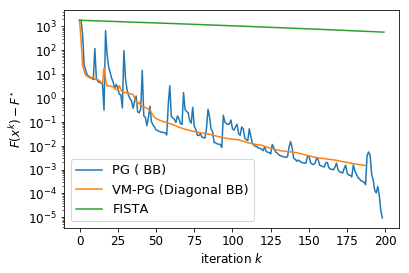}
  \caption{$\ell_1$ penalized LS for CIFAR using
  ($N,n$)=($50000,3072$)}
  \label{fig:sub5}
\end{subfigure} \hspace{0.05\textwidth} %
\begin{subfigure}{.4\textwidth}
  \centering
  \includegraphics[width=1\linewidth]{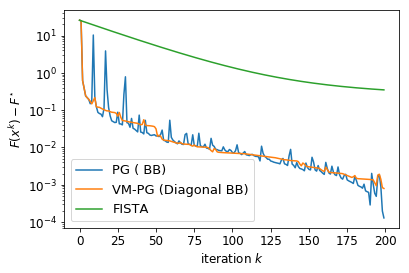}
  \caption{$\ell_1$ penalized LS for CIFAR using
   ($N,n$)=($625,3072$)}
  \label{fig:sub6}
\end{subfigure}
\vspace{0mm}
\caption{Additional convergence comparison between VM-PG with diagonal BB stepsize, PG with BB stepsize, and Accelerated PG (FISTA) for penalized least square (LS), penalized logistic regression (LR). }
\vspace{-5pt}
\label{fig:test}
\end{figure*}

\subsection{Additional results for the real-world datasets} 
Fig. \ref{fig:test} shows that VM-PG (DBB) is not worse and often faster than PG (BB) exhibiting stable behavior (less oscillations) to convergence. We also provide the results for FISTA as a state-of-art baseline. For the LS problems FISTA is much slower  than PG (BB) and VM-PG (DBB); however, but is often faster for the LR problem. This however, can be remedied through careful selection of the $\mu$ parameter in the VMPG algorithm. Such, optimal selection of the $\mu$ parameter for improved performance of VMPG compared to FISTA is an on-going research topic. 



\end{document}